\documentclass[a4paper]{article}

\usepackage{delarray,verbatim,enumerate}
\usepackage{amsmath,amsthm,amstext,amsbsy,amssymb,amsfonts,amscd}
\usepackage{MnSymbol}
\usepackage[latin1]{inputenc}
\usepackage[francais]{babel}
\usepackage[T1]{fontenc}

\usepackage{color}
\usepackage[colorlinks,linkcolor=blue,urlcolor=blue,linktocpage,dvips]{hyperref}

\usepackage[all]{xy}%

\newtheorem{Th}{Théorème}%[section]
\newtheorem{Lem}[Th]{Lemme}
\newtheorem{Prop}[Th]{Proposition}
\newtheorem{Cor}[Th]{Corollaire}

\newtheorem{Sco}[Th]{Scolie}
\newtheorem{Def} [Th]{Définition}

\newcommand{\etoile}{\mbox{\fontsize{15}{24}\selectfont $*$}}

\def\Preuve{\smallskip\noindent {\it Preuve.~}}
\def\PreuveTh{\smallskip\noindent {\it Preuve du Théorème.~}}
\def\Remarque{\smallskip\noindent {\it Remarque.~}}

\def\ie{{\it i.e. }}	\def\cf{{\it cf. }}	\def\eg{{\it e.g. }}
	
	\def\N{\mathbb N}
\def\Z{\mathbb {Z}}     			\def\Q{\mathbb Q}
\def\Zl{\mathbb{Z}_\ell}			

\def\L{\Lambda}					

\def\L{\mathcal  L}	      			\def\T{\mathcal  T}
\def\J{\mathcal  J}  		\def\C{\mathcal  C}		\def\R{\mathcal  R}
\def\E{\mathcal  E} 		\def\G{\mathcal  G}		\def\U{\mathcal  U}
\def\G{\mathcal  G}				
 			  	\def\Cl{\mathcal  C\ell}

\def\p{\mathfrak p}		\def\q{\mathfrak q}	\def\l{\mathfrak l}		\def\r{\mathfrak r}

\def\wi{\widetilde}		

\def\Gal{\operatorname{Gal}}

\date{}

\title{\huge Propagation de la 2-birationalité}

\author{{\small par}\\
\\
Claire {\sc Bourbon} \& Jean-François {\sc Jaulent}}

\begin{document}

\maketitle
\bigskip

{\footnotesize
\noindent{\bf Résumé.} Nous étudions la propagation de la 2-birationalité dans les 2-extensions de corps de nombres. Nous prouvons que pour toute extension quadratique totalement imaginaire 2-birationnelle  $L$ d'un corps de nombres 2-rationnel totalement réel $K$, la propagation de la 2-birationalité par 2-extension de $K$ n'est possible, à composition près par une 2-extension cyclotomique, que dans le cas quadratique. Nous la caractérisons complètement en termes de ramification modérée ce qui permet de construire des tours infinies de telles 2-extensions.}
\medskip

{\footnotesize
\noindent{\bf Abstract.}  Let $L/K$ be a 2-birational CM-extension of a totally real 2-rational number field. We characterize in terms of tame ramification totally real 2-extensions $K'/K$ such that the compositum $L'=LK'$ is still 2-birational. In case the 2-extensions $K'/K$ is linearly disjoint from the cyclotomic $\Z_2$-extension $K^c/K$, we prove that $K'/K$ is at most quadratic. In the other hand we construct infinite towers of such 2-extensions.}
\bigskip

\section{\large Introduction}

La notion de corps $S$-rationnel a été introduite dans \cite{JS1}, en liaison 
avec les résultats de \cite{Wi}, pour généraliser la notion de corps de 
nombres $\ell$-rationnel rencontrée implicitement dans des contextes 
variés par plusieurs auteurs puis explicitement définie et étudiée
 par \cite{MN} d'une part et  \cite{GJ} d'autre part (\cf  \cite{JN}). Rappelons ce dont il 
 s'agit~: si $\ell$ est un nombre premier et $S$ un sous-ensemble non vide 
 de l'ensemble $Pl^{\;\ell}_K$ des places de $K$ au-dessus de $\ell$, on dit que le corps de
 nombres $K$ est $S$-rationnel   lorsque le groupe de Galois $\G_K=\Gal(M'/K)$ 
 de sa pro-$\ell$-extension (galoisienne) $\ell$-ramifiée maximale est le 
pro-$\ell$-produit libre 
$$
\G_K \simeq (\ \etoile\hspace{-12pt}
\underset{\underset{{\p}\not\in S}{{\p}|\ell\infty} {\ }} \ {\G}_{K_\p} ) \etoile {\mathcal F}\leqno{(i)}
$$

des groupes de Galois locaux $\G_ {K_\p} =Gal(\bar{ K}_ \p /K_ \p)$ 
respectivement attachés aux pro-$\ell$-extensions maximales $\bar K_ {\mathfrak p} $ 
des completés $K_ {\mathfrak p} $ de $K$ aux places réelles ou $\ell$-adiques qui 
ne sont pas dans $S$, et d'un pro-$\ell$-groupe libre $\mathcal F$ ; dans ce cas,
le nombre $f$ de générateurs de $\mathcal F$ est donné par la formule 
$$f=d-r-c-l+s+1,\leqno{(ii)}$$
o\`u $d$ est la somme des degrés locaux 
$d=\sum_{{\mathfrak l}\in S}[K_{{\mathfrak l}} :{\mathbb Q}_\ell]$ 
et $r$, $c$, $l$, $s$ sont respectivement les nombres de places
 réelles , complexes, $\ell$-adiques ou dans $S$ de $K$ (cf. \cite{JS1}), th 2.7). 
 Lorsque 
$S$ est un singleton $\{{\mathfrak l} \}$, on parle de corps ${\mathfrak l}$-rationnel et 
si  $Pl^{\;\ell}_K$ est lui même un singleton, on dit tout simplement que $K$ est
$\ell$-rationnel.\smallskip

Dans ce dernier cas (\ie pour $l=s=1$), la formule $(ii)$ ci-dessus donne $f=c+1$; et la condition $(i)$ affirme que $\G_K$ est le pro-$\ell$-produit libre d'un pro-$\ell$-groupe libre de dimension $c+1$ (qui s'identifie au groupe de Galois de la sous-extension $\infty$-décomposée maximale $M/K$ de $M'/K$) et des $r$ sous-groupes de décomposition attachés aux places réelles de $K$. Lorsque le corps $K$ contient en outre les racines $\ell$-ièmes de l'unité, ceci a lieu si et seulement si le $\ell$-groupe $C\ell'_K$ des $\ell$-classes de diviseurs de $K$ (\ie le quotient du $\ell$-groupe des classes d'idéaux prises au sens restreint par le sous-groupe  engendré par la classe de l'idéal premier au dessus de $\ell$) est trivial  (\cf \eg \cite{JN}), ce qui s'écrit:
$$\Cl'_K=1,\leqno{(iii)}$$
de sorte que la notion de $\ell$-rationalité coïncide alors avec celle de $\ell$-régularité introduite par Kummer dans l'étude des corps cyclotomiques ${\mathbb Q}[\zeta_\ell]$.\medskip

La question de la propagation de la $S$-rationalité dans une $\ell$-extension $L/K$ de corps de nombres a été complètement résolue dans \cite{JS1} pour les $\ell$ premiers impairs. Pour $\ell=2$ et $L/K$ quadratique il peut arriver que le corps de base $K$ soit $\mathfrak{l}$-rationnel en une place $2$-adique décomposée dans l'extension $L/K$ et que $L$ soit $\mathfrak{l}$-rationnel en chacune des deux places au-dessus de $\mathfrak{l}$; on dit alors que le corps $L$ est {\em birationnel}. Le passage de la $2$-rationalité à la $2$-birationalité dans une 2-extension $L/K$ a été complètement traité dans \cite{JS2} pour $K$ totalement réel. Rappelons le résultat principal: \smallskip

\setcounter{Th}{-1}
\begin{Th}
Soit $L/K$ une extension quadratique totalement imaginaire d'un corps de nombres totalement réel. Les assertions suivantes sont équivalentes : 
\begin{enumerate}
\item[(i)] Le corps $L$ est $2$-birationnel,
\item[(ii)] Le corps $K$ est $2$-rationnel, son unique place $2$-adique est décomposée dans $L/K$ et l'extension $L/K$ est ramifiée modérement soit en une place semi-primitive $\mathfrak{p}$ soit en deux places primitives $\mathfrak{p}$ et $\mathfrak{q}$.  
\end{enumerate}
\end{Th}

Dans ce contexte, une place finie $\p$ du corps de base $K$ est dite {\em primitive} lorsqu'elle est totalement inerte dans la $\Z_2$-extension cyclotomique $K^c/K$; {\em semi-primitive} lorsqu'elle est décomposée dans le premier étage de $K^c/K$ et inerte au-delà.  Une telle place est, en particulier, modérée (\ie ici étrangère à 2).
\medskip

L'objet du présent travail est d'étudier la propagation de la 2-birationalité par 2-extension (galoisienne) du corps de base. Le problème est le suivant: étant donnés une extension quadratique à conjugaison complexe $L$ d'un corps totalement réel $K$ d'une part, une $2$-extension totalement réelle $K'$ de $K$ d'autre part, nous regardons à quelle condition sur les extensions $L/K$ et $K'/K$ la $2$-birationalité de $L/K$ se propage à l'extension induite $LK'/K'$.\smallskip

Le résultat {\em a priori } surprenant que nous obtenons est le suivant:

\begin{Th}\label{ThPr}
Lorsque la 2-extension totalement réelle $K'/K$ est ramifiée modérément en une place $\p$ (auquel cas celle-ci est primitive et c'est l'unique place modérée qui se ramifie dans  $K'/K$), la propagation de la 2-birationalité de $L/K$ à $LK'/K'$ a lieu si et seulement si  les deux conditions qui suivent sont réalisées:
	\begin{enumerate}
	\item[(i)] l'extension $L/K$ est ramifiée modérément en exactement deux places: en la place $\mathfrak{p}$ et en une autre place primitive $\mathfrak{q}$; 
	\item[(ii)] le corps $K'$ provient, par composition avec un étage fini $K_n$ de  la $\Z_2$-extension cyclotomique $K^c/K$, d'une extension quadratique $K''$ de $K$.
	\end{enumerate}  
\end{Th}

En d'autres termes, si l'on impose à $K'/K$ d'être linéairement disjointe de la $\Z_2$-extension cyclotomique $K^c/K$, la propagation de la birationalité est impossible lorsque $L/K$ est modérément ramifiée en une place semi-primitive;  et lorsque $L/K$ est modérément ramifiée en deux places primitives $\p$ et $\q$, elle n'est possible que par extension quadratique $\p$-modérément ou $\q$-modérément ramifiée, ce qui conduit, comme nous le verrons, à deux possibilités (et deux seulement) pour $K'$. En itérant le processus, en revanche, il est alors facile de construire des tours infinies de telles extensions.

%%%%%%%%%%%%%%%%%%%%%%%%%%%%%%%%%%%%%%%%%%%%%%%%%%%%%%%
\section{\large Corps 2-rationnels}
%%%%%%%%%%%%%%%%%%%%%%%%%%%%%%%%%%%%%%%%%%%%%%%%%%%%%%%

Pour la commodité du lecteur nous rappelons brièvement ci-dessous quelques uns des résultats de \cite{GJ} et \cite{MN} sur les notions de régularité et de rationalité:

\begin{Def}
Soit $\ell$ un nombre premier. Un corps de nombres $K$ est dit:
\begin{enumerate}
\item[(i)] $\ell$-régulier, lorsque le $\ell$-sous-groupe de Sylow du noyau $R_2(K)$ dans $K_2(K)$ des symboles réguliers est trivial;
\item[(ii)] $\ell$-rationnel, si le groupe de Galois $\,\G_K=\Gal(M/K)$ de sa pro-extension maximale  $\ell$-ramifiée $\infty$-décomposée est un pro-$\ell$-groupe libre.
\end{enumerate}
\end{Def}

Lorsque $K$ contient le sous-corps réel maximal $\Q[\zeta \,+\,\zeta^{-1}]$ du $\ell$-ième corps cyclotomique $\Q[\zeta]$, il est équivalent d'affirmer qu'il est $\ell$-régulier ou qu'il est $\ell$-rationnel (\cf \cite{JN}, Théorème 1.2). C'est évidemment le cas pour $\ell=2$. Ainsi:

\begin{Th}\label{Car2r}
$K$ est 2-rationnel lorsqu'il vérifie les propriétés équivalentes:
\begin{enumerate}
\item[(i)] La 2-partie du noyau dans $K_2(K)$ des symboles réguliers est trivial.
\item[(ii)] Le groupe de Galois  $\G_K=\Gal(M/K)$ de sa pro-2-extension maximale de $K$ 2-ramifiée et $\infty$-décomposée est un pro-2-groupe libre (nécessairement sur 1+$c_K$ générateurs, si $c_K$ est le nombre de places complexes de $K$).
\item[(iii)] Le groupe de Galois  $\G^{ab}_K=\Gal(M^{ab}/K)$ de sa pro-$\ell$-extension abélienne maximale 2-ramifiée $\infty$-décomposée est un $\Z_2$-module libre de rang 1+$c_K$.
\item[(iv)] Le corps $K$ vérifie la conjecture de Leopoldt (pour le nombre premier 2) et le sous-module de torsion $\T_K$ de $\Gal(M^{ab}/K)$ est trivial.
\item[(v)]  $K$ possède une unique place dyadique $\l$ et son 2-groupe $\,\Cl'_K$ des 2-classes d'idéaux au sens restreint est trivial.
\end{enumerate}
\end{Th}
\smallskip

L'équivalence des diverses assertions {\em (i)} à {\em (v)} n'est autre que la déclinaison pour $\ell=2$ du Th. 1.2 de \cite{JN}. La formulation {\em (ii)} montre clairement que la 2-rationalité se propage par 2-extension 2-ramifiée $\infty$-décomposée, donc à chaque étage fini $K_n$ de la $\Z_2$-extension cyclotomique. En particulier les 2-groupes de 2-classes $\,\Cl'_{K_n}$  sont tous triviaux et le 2-groupe $\,\wi\Cl_K$ des 2-classes logarithmiques de $K$ est ainsi trivial. Comme expliqué dans \cite{Ja2}, il suit de là qu'un tel corps satisfait aussi la conjecture de Gross (pour le nombre premier 2).\smallskip

La question de la propagation de la $\ell$-rationalité a été complétement résolue par \cite{GJ} et \cite{MN}. Elle s'appuie sur la notion de place primitive. Pour $\ell=2$, on a:

\begin{Def}
Étant donné un corps de nombres $K$, un ensemble $S$ de places modérées de $K$ (\ie ici de places finies étrangères à 2) est dit primitif (relativement au premier 2) lorsque la famille des logarithmes de Gras de ces places (\ie de leurs images dans le groupe de Galois $\L=\Gal (Z/K)$ du compositum $Z$ des $\Z_2$-extensions de $K$) peut être complétée en une $\Z_2$-base de $\L$.
\end{Def}  

Ainsi, lorsque $K$ un corps totalement réel qui satisfait la conjecture de Leopoldt (pour le premier 2), par exemple un corps 2-rationnel, les ensembles primitifs de places modérées de $K$ sont exactement les singletons $S=\{\l\}$, où $\l$ est une place de $K$ totalement inerte dans la $\Z_2$-extension cyclotomique $Z/K$.

\begin{Def}
Une 2-extension $L/K$ est dite primitivement ramifiée lorsque l'ensemble $\mathcal{R}_{L/K}$ des places modérément ramifiées dans $L/K$ est primitif.
\end{Def}

Le résultat de propagation (\cf \cite{JN}, Théorème 3.5) s'énonce alors comme suit:

\begin{Th}\label{Propa2r}
Étant donné une 2-extension de corps de nombres $L/K$, les deux conditions suivantes sont équivalentes : 
\begin{enumerate}
\item[(i)] Le corps $L$ est 2-rationnel.
\item[(ii)] Le corps $K$ est 2-rationnel et l'extension $L/K$ est primitivement ramifiée.
\end{enumerate}
\end{Th}

Donnons par exemple la liste des corps multiquadratiques 2-rationnels:

\begin{Cor}\label{cmq2r}
Les corps multiquadratiques réels $K$ qui sont 2-rationnels sont les sous-corps des corps biquadratiques $\Q[\sqrt{2}$,$\sqrt{p}]$ pour les nombres premiers $p$ primitifs, c'est-à-dire vérifiant $p\equiv\pm3 \mod8$. En d'autres termes ce sont: 
\begin{itemize}
\item[(i)] le corps des rationnels $\Q$; 
\item[(ii)]  les corps quadratiques $\Q[\sqrt{2}]$, $\Q[\sqrt{p}]$, $\Q[\sqrt{2p}]$;
\item[(iii)]  les corps biquadratiques $\Q[\sqrt{2},\sqrt{p}]$;
\end{itemize}
où $p$ est un premier primitif: $p\equiv\pm3 \mod8$.\smallskip

\noindent Les corps multiquadratiques imaginaires $2$-rationnels sont:
\begin{itemize} 
\item[(iv)]  les corps quadratiques $\Q[\sqrt{-1}]$, $\Q[\sqrt{-2}]$, $\Q[\sqrt{-p}]$, $\Q[\sqrt{-2p}]$;
\item[(v)] les corps biquadratiques $\Q[\sqrt{-1},\sqrt{2}]$, $\Q[\sqrt{-1},\sqrt{p}]$, $\Q[\sqrt{-2},\sqrt{p}]$;
\item[(vi)] les corps triquadratiques $\Q[\sqrt{-1},\sqrt{2},\sqrt{p}]$;
\end{itemize}
où $p$ est un premier primitif: $p\equiv\pm3 \mod8$.
\end{Cor}

\Preuve Le corps $\Q\,$ étant 2-rationnel, il suffit d'écrire que la ramification modérée ne peut concerner qu'au plus un premier impair $p$, lequel doit en outre être primitif, \ie inerte dans le premier étage $\Q[\sqrt{2}]/\Q$ de la $\Z_2$-extension de $\Q$; ce qui se traduit par la congruence annoncée: $p\equiv\pm3\mod8$.

%%%%%%%%%%%%%%%%%%%%%%%%%%%%%%%%%%%%%%%%%%%%%%%%%%%%%%%
\section{\large Corps 2-birationnels}
%%%%%%%%%%%%%%%%%%%%%%%%%%%%%%%%%%%%%%%%%%%%%%%%%%%%%%%

Venons-en maintenant à la notion de corps 2-birationnel (\cf \cite{JS1} et \cite{JS2}):

\begin{Def}
Un corps de nombres $K$ est dit $\l$-rationnel en une place $\l$ au-dessus de $\ell$ si sa $\ell$-extension  $\ell$-ramifiée, $\l$-décomposée maximale est triviale.
\end{Def}

Naturellement la trivialité du groupe de Galois de la pro-$\ell$-extension  $\ell$-ramifiée, $\l$-décomposée maximale se lit sur son abélianisé. On peut alors interpréter la notion de la pro-$\l$-rationalité en termes de corps de classes. Introduisons pour cela quelques notations de la Théorie $\ell$-adique du corps de classes (\cf \cite{Ja2}).\smallskip

Notons:\smallskip

\begin{itemize}
\item $\R_{K_\p}$ le compactifié $\ell$-adique du groupe multiplicatif du complété $K_\p$;
\item $\mu_\p$ le sous-groupe de torsion de $\R_{K_\p}$;
\item $\J_K=\prod_\p^{res}\R_{K_\p}$ le $\ell$-adifié du groupe des idèles de $K$;
\item $\R_K=\Zl\otimes K^\times$ son sous-groupe principal ;
\item $\C_K=\J_K/\R_K$ enfin le $\ell$-groupe des classes d'idèles.
\end{itemize}\smallskip

Le groupe $\C_K$ s'identifie au groupe de Galois de la pro-$\ell$-extension abélienne maximale de $K$. Et d'après  \cite{JS1} Th. 1.7, il vient:

\begin{Prop}
Le corps $K$ est $\l$-rationnel si et seulement si on a l'identité: \smallskip

\centerline{$ \J_K = \R_K \R_{K_\l}\prod\limits_{\p\nmid{\ell\infty}}\mu_\p$;}\smallskip

\noindent ce qui a lieu si et seulement si les deux conditions suivantes se trouvent réunies : 
\begin{enumerate}
\item[(i)] le groupe des $\ell$-classes d'idéaux $\,\Cl_K$ du corps $K$ est trivial;
\item[ii)] l'application de semi-localisation $s'_\l$ induit une surjection du tensorisé $\E_K'$ du groupe des $\ell$-unités de $K$ sur le produit $\R_l'=\prod_{\p\mid{\ell\infty};\p\ne\l}\R_{K_\p}$.% étendu aux places $\mathfrak{l}'\neq\mathfrak{l}$ divisant $\ell$ et aux places réelles.
\end{enumerate}
\end{Prop}

\Remarque
Sous les conditions précédentes, $K$ est logarithiquement principal en ce sens que son $\ell$-groupe des classes logarithiques $\wi\Cl_K$ est trivial.  En particulier $K$ vérifie banalement la conjecture de Gross généralisée (pour le premier $\ell$).
\medskip

Ces définitions générales étant posées, nous pouvons spécifier au cas $\ell=2$.

\begin{Th}
Un corps totalement imaginaire $L$ est dit $2$-birationnel lorsqu'il est $\mathfrak{l}$-rationnel  en chacune des places $2$-adiques (en ce sens qu'il n'admet pas de $2$-extension, $2$-ramifiée, $\mathfrak{l}$-décomposée non triviale), ce qui a lieu si et seulement si les trois conditions suivantes sont vérifiées :
\begin{itemize}
\item[(a)] $L$ admet exactement $2$ places $2$-adiques $\mathfrak{l}$ et $\mathfrak{l'}$; 
\item[(b)] le $2$-groupe $\Cl_L$ des $2$-classes d'idéaux de $L$ est trivial, \ie $2$-groupe des classes d'idéaux de $L$ est engendré par les images de $\l$ et de $\l'$.
\item[(c)] les plongements canoniques de $L^\times$ dans $L^\times_\l$ et $L^\times_{\l'}$ induisent des isomorphismes $\E'_L\simeq \R_{L_\l}\simeq \R_{L_{l'}}$ du tensorisé $2$-adique $\E'_L$ du groupe des $2$-unités de $L$ sur les compactifiés locaux associés aux places $2$-adiques.
\end{itemize}
\end{Th}

Nous avons maintenant besoin d'introduire la notion de place semi-primitive. 

\begin{Def}
Soit $K$ un corps de nombres totalement réel qui satisfait la conjecture de Leopoldt en $\ell=2$ et $K^c$ sa $\Z_2$-extension cyclotomique. Nous disons qu'une place finie modérée (\ie ici ne divisant pas 2) est:
\begin{itemize}
\item {\em primitive}, lorsque son image dans le groupe procyclique $Gal(K^c/K)$ n'est pas un carré, autrement dit lorsqu'elle n'est pas décomposée dans $K^c/K$ ; 
\item {\em semi-primitive}, lorsqu'elle se décompose dans le premier étage de la $\Z_2$-extension cyclotomique $K^c/K$ mais pas au-delà.
\end{itemize}
\end{Def}

Cela posé, nous avons les résultats suivants  (cf \cite{JS2} Th. 4 \& Prop. 5):

\begin{Th}\label{Carac2bi}
Soit $L/K$ une extension quadratique totalement imaginaire d'un corps de nombres totalement réel. Les assertions suivantes sont équivalentes : 
\begin{enumerate}
\item[(i)] Le corps $L$ est $2$-birationnel.
\item[(ii)] Le corps $K$ est $2$-rationnel, son unique place $2$-adique est décomposée dans $L/K$ et l'extension $L/K$ est ramifiée modérement soit en une place semi-primitive $\mathfrak{p}$ soit en deux places primitives $\mathfrak{p}$ et $\mathfrak{q}$.  
\item[(iii)] $K$ est $2$-rationnel; $L$ est $2$-logarithmiquement principal; et $L/K$ est $2$-dé\-composée et ramifiée modérément en au moins une place.
\end{enumerate}
\end{Th}

%%%%%%%%%%%%%%%%%%%%%%%%%%%%%%%%%%%%%%%%%%%%%%%%%%%%%%%
\section{\large Corps multiquadratiques 2-birationnels}
%%%%%%%%%%%%%%%%%%%%%%%%%%%%%%%%%%%%%%%%%%%%%%%%%%%%%%%
\
Nous nous nous proposons dans cette section de mettre en \oe uvre les équivalences données par le Th. \ref{Carac2bi} et la classification des corps multiquadratiques réels 2-rationnels explicitée dans le Cor. \ref{cmq2r} pour dresser la liste des corps multiquadratiques 2-birationnels: un tel corps $L$ s'obtient, en effet, par composition d'un corps quadratique imaginaire $k=\Q[\sqrt{-d}]$ et d'un corps multiquadratique réel $K$, qui est nécessairement 2-rationnel en vertu de la condition $(ii)$ du Théorème.\smallskip

Le résultat que nous obtenons est le suivant:

\begin{Th}\label{cmq2bir}
Soient $L$ un corps multiquadratique imaginaire et $K$ son sous-corps réel maximal. Le corps $L$ est alors le compositum du corps  multiquadratique réel $K$ et d'un corps quadratique imaginaire $k=\Q[\sqrt{-d}]$, pour un $d\geq 1$ sans facteur carré. Et $L$ est 2-birationnel si et seulement si $K$ est 2-rationnel (\ie un sous-corps de $\Q[\sqrt{2},\sqrt{p}]$ avec  $p\equiv \pm3\mod8$) et que l'on est dans l'une des quatre configurations suivantes:
\begin{enumerate}
\item[(a)] Pour $K[\sqrt{2}]=\Q[\sqrt{2}]$: 
	\begin{enumerate}
	\item[(i)]$d=q$ premier avec $q\equiv7\mod16$;
	\item[(ii)]$d=qq'$ avec $q,q'$ premiers et $q\equiv-q'\equiv\pm3\mod8$.
	\end{enumerate}
\item[(b)] Pour $K[\sqrt{2}]=\Q[\sqrt{2},\sqrt{p}]$:
	\begin{enumerate}
	\item[(i)]  $d=q\neq{p}$ premier avec $-q\equiv p \equiv \pm 3 \mod8$ et $(\frac{p}{q})=-1$;
	\item[(ii)] $d=q\neq{p}$ premier avec $-q\equiv p \equiv \pm 3 \mod8$ et $(\frac{p}{q})=+1$.
	\end{enumerate}
\end{enumerate}
\end{Th}

Pour démontrer cela, nous nous appuierons sur le lemme:

\begin{Lem}
\label{equiv}
Soit $L$ une extension quadratique totalement imaginaire d'un corps totalement réel $K$. Pour tout $n\in\N$, soit $K_n$  le $n$-ième étage de la $\Z_2$-extension cyclotomique $K^c$ de $K$ et $L_n=K_nL$ le compositum. On a alors les équivalences :
\item [(i)] $K$ est $2$-rationnel $\Longleftrightarrow$ $K_n$ est $2$-rationnel.
\item [(ii)] $L$ est $2$-birationnel $\Longleftrightarrow$ $L_n$ est $2$-birationnel.
\end{Lem}

\Preuve Commençons par établir le Lemme \ref{equiv}.

 Observons d'abord que la $\Z_2$-extension cyclotomique $K^c/K$ étant 2-ramifiée, il en est de même de toutes ses sous-extensions finies $K_n/K$, de sorte que la 2-rationalité de $K$ est bien équivalente à celle de chacun des $K_n$ en vertu du Théorème \ref{Propa2r}; d'où l'équivalence $(i)$ du Lemme.
 
 Ce point acquis, intéressons-nous à la 2-birationalité de $L$. D'après la formulation $(iii)$ donnée par le Théorème \ref{Carac2bi}, celle-ci se caractérise par la 2-rationalité de $K$, la 2-décomposition de $L/K$ et le trivialité du 2-groupe des classes logarithmiques $\,\wi\Cl_L$; et il s'agit simplement de vérifier que ces conditions se propagent le long de la tour cyclotomique $L^c/L$. Or:
 \begin{itemize}
 \item[(i)] La 2-rationalité de $K_n$ se lit indifféremment à n'importe quel étage de la tour en vertu du point $(i)$ ci-dessus.
  \item[(ii)] Si  $K$ est 2-rationnel, il admet une unique place dyadique $\l$ et il en est de même pour chacun des $K_n$;  l'unique place dyadique $\l_n$ de $K_n$ se décompose dans $L_n/K_n$ si et seulement si $\l$ est décomposée dans $L/K$.
 \item[(iii)] Enfin, le 2-groupe des classes logarithmiques $\wi\Cl_L$ étant le quotient des copoints fixes par $\Gamma=\Gal(L^c/L)$ de l'abélianisé $\C'_L$ du groupe de Galois de la pro-2-extension totalement décomposée maximale $\wi L$ de $L$ (\cf \cite{Ja1, Ja2}), la trivialité de $\wi\Cl_L$ affirme tout simplement l'égalité $\wi L= L$; ce qui se lit encore à n'importe quel étage de la tour, puisque $\wi L$ est aussi la pro-2-extension totalement décomposée maximale de chacun des $L_n$.
 \end{itemize}
 
 \PreuveTh 
Si $L$ est 2-birationnel, son sous-corps réel $K$ est 2-rationnel \ie (par le Corollaire \ref{cmq2r}) contenu dans un corps biquadratique $\Q[\sqrt2,\sqrt p]$ pour un premier primitif $p\equiv \pm 3 \mod 8$. Cela étant, le Lemme \ref{equiv} nous permet de remplacer $K$ par n'importe quel étage fini de sa $\Z_2$-extension cyclotomique, donc de supposer par exemple $K\supset \Q[\sqrt 2]$. Distinguons deux cas:\medskip

$\bullet$ Pour $K[\sqrt{2}]=\Q[\sqrt{2}]$: \smallskip

\noindent Dans ce cas, toujours d'après le Lemme \ref{equiv}, nous pouvons remplacer $K$ par $K[\sqrt 2]=\Q[\sqrt 2]$ donc par $\Q$; de sorte que le corps $L$ est 2-birationnel si et seulement si $k=\Q[\sqrt{-d}]$ l'est; ce qui suppose, d'après le Théorème \ref{Carac2bi} $(ii)$:\smallskip

\begin{itemize}
\item[(i)] 2 décomposé dans l'extension $k/\Q$, \ie $d\equiv -1\mod 8$, et;\smallskip

\item[(ii)] $k/\Q$ ramifiée modérément soit en une place semi-primitive $q$, \ie $d=q$ ou $2q$, avec $q\equiv \pm 7 \mod 16$; soit en deux places primitives $q$ et $q'$, \ie $d=qq'$ ou $2qq'$, avec $q \equiv \pm q' \equiv \pm 3 \mod 8$.
\end{itemize}\smallskip

\noindent En fin de compte, il vient $d=q \equiv 7 \mod 16$ ou $d=qq'$ avec $q\equiv -q' \equiv \pm 3 \mod 8$.\medskip

$\bullet$ Pour $K[\sqrt{2}]=\Q[\sqrt 2, \sqrt p]$: \smallskip

\noindent Dans ce cas, toujours d'après le Lemme \ref{equiv}, nous pouvons remplacer $K$ par $\Q[\sqrt p]$ et $L$ par $\Q[\sqrt{p},\sqrt{-d}]$; d'où par le Théorème \ref{Carac2bi} $(ii)$, les trois conditions:\smallskip

\begin{itemize}
\item[(i)] Le corps quadratique $K=\Q[\sqrt p]$ est 2-rationnel, \ie on a: $p\equiv \pm 3\mod 8$.\smallskip

\item[(ii)]  Son unique place dyadique $\l$ est décomposée dans $L/K$;  autrement dit soit 2 est décomposé dans l'extension $k/\Q$, auquel cas on a: $d\equiv -1\mod 8$; soit 2 est décomposé dans l'extension $\Q[\sqrt{-dp}]/\Q$, auquel cas on a: $dp \equiv -1 \mod 8$, conformément au schéma d'extensions:
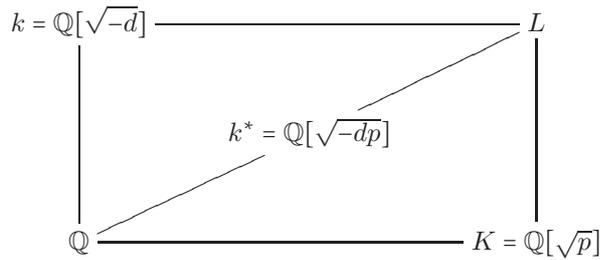
\begin{displaymath}
 \xymatrix { 
    k=\Q[\sqrt{-d}] \ar@{-}[rr] \ar@{-}[dd] && L\ar@{-}[dd] \\
    & k^*=\Q[\sqrt{-dp}] \ar@{-}[ur] \\
    \Q \ar@{-}[rr] \ar@{-}[ur] && K=\Q[\sqrt{p}]
    }
\end{displaymath}
Et, dans les deux éventualités, $d$ est donc impair.\smallskip

\item[(iii)] $L/K$ se ramifie modérément soit en une unique place semi-primitive $\q$, soit en deux places primitives $\q$ et $\q'$.
\begin{itemize}
\item Dans le premier cas la place $\q$ provient alors d'un premier $q\neq p$ ramifié dans $k/\Q$  mais inerte dans $K/\Q$, auquel cas on a: $d=q$ ou $d=pq$ et $(\frac{p}{q})$=$-1$. Reste simplement à vérifier la semi-primitivité de $\q$. Le schéma

\begin{displaymath}
 \xymatrix { 
    &K \ar@{-}[d] ^{in}\ar@{-}[r] ^{dec}&K[\sqrt 2]\ar@{-}[d]^{dec}\ar@{-}[r]^{in}  \ar@{-} &K[\sqrt {2+\sqrt 2}]\ar@{-}[d]^{dec}\ar@{.}[r]&K^c \ar@{.}[d]\\
    &\Q \ar@{-}[r] ^{in}&\Q[\sqrt 2] \ar@{-}[r]^{in}&\Q[\sqrt {2+\sqrt 2}]\ar@{.}[r]&\Q^c\\
    }
\end{displaymath}
(où est indiqué le comportement des places au-dessus de $q$) montre alors que celle-ci correspond à la primitivité de $q$, qui s'écrit: $q\equiv\pm3\mod 8$.\smallskip

\item Dans le second cas, $\q$ et $\q'$, du fait de leur primitivité, proviennent nécessairement d'un même premier $q\ne p$ ramifié dans $k/\Q$ et décomposé dans $K/\Q$, auquel cas on a: $d=q$ ou $d=pq$ et $(\frac{p}{q})$=$+1$, ainsi que la congruence $q\equiv\pm 3$ qui traduit la primitivité de $q$.
\end{itemize}\smallskip

\noindent En fin de compte, quitte à échanger $k$ et $k^*$, on obtient $d=q\neq{p}$ premier avec $-q\equiv p \equiv \pm 3 \mod8$ et $(\frac{p}{q})=-1$ dans le premier cas; $(\frac{p}{q})=+1$ dans le second; ce qui achève la démonstration.
\end{itemize}\medskip

\Remarque Le cas $(a)$ du Théorème \ref{cmq2bir} redonne naturellement la classification des corps quadratiques imaginaires 2-birationnels donnée dans \cite{JS1}  (Cor.  1.12), ainsi que la liste des corps quadratiques imaginaires 2-logarithmiquement principaux établie dans \cite{So} (restreinte ici à ceux qui sont 2-décomposés):
     
\begin{Cor}\label{cq2bir}
Les corps quadratiques $2$-birationnels sont les corps quadratiques imaginaires $K=\Q[\sqrt{-p}]$ pour $p$ premier de la forme $p\equiv7\mod16$, et $K=\Q[\sqrt{-pq}]$ pour $p$ et $q$ premiers $p\equiv-q\equiv3\mod 8$.
\end{Cor}

%%%%%%%%%%%%%%%%%%%%%%%%%%%%%%%%%%%%%%%%%%%%%%%%%%%%%%%%
\section{\large Propagation de la birationalité}
%%%%%%%%%%%%%%%%%%%%%%%%%%%%%%%%%%%%%%%%%%%%%%%%%%%%%%%%

Dans le cas $(b)$ du Théorème \ref{cmq2bir}, on a $d=q$ ou $d=-pq$, de sorte que le corps $L=K[\sqrt{-d}]$ provient, par composition avec $K$ , indifféremment du corps quadratique imaginaire $k=\Q[\sqrt{-q}]$ comme du corps  $k^*=\Q[\sqrt{-pq}]$. Il est intéressant d'observer que, du fait des congruences satisfaites par $p$ et $q$, un et un seul d'entre eux (à savoir $k^*$) se trouve être birationnel. Plus précisément, $k^*$ est  une extension 2-birationnelle de $\Q$ qui est ramifiée modérément en deux places primitives, dont l'une se ramifie dans $K/\Q$ tandis que l'autre se décompose.\medskip

Nous allons voir que cette sitation est, en fait, générale:

\begin{Th} \label{ThPropa}
Soit $K$ un corps 2-rationnel totalement réel; $L=K[\sqrt{-\delta}]$ (avec $\delta>>0$ dans $K$) une extension quadratique totalement imaginaire et $K'$ une 2-extension totalement réelle non triviale de $K$ linéairement disjointe de la $\Z_2$-extension cyclotomique $K^c$; soit enfin $L'=LK'$ le compositum.de $L$ et de $K'$. Alors:
\begin{itemize}
\item [(i)] Si $L'/K'$ est 2-birationnelle, l'extension $L/K$ ne peut elle-même être 2-birationnelle que si sont réunies les deux conditions suivantes:
\begin{itemize}
\item [(a)]  L'extension $K'/K$ est aussi quadratique: $[K':K]=2$.
\item [(b)]  L'extension  $L/K$ vérifie l'une des deux propriétés qui suivent:
\begin{itemize}
\item [(b1)]  ou bien $L/K$ est ramifiée modérément en une unique place primitive $\q$, laquelle est inerte dans $K'/K$;
\item [(b2)]  ou bien $L/K$ est  ramifiée modérément en exactement deux places primitives $\p$et $\q$ et l'extension quadratique $K'/K$ est ramifiée modérément en exactement l'une de ces deux places et  décomposée en l'autre.
\end{itemize}
\end{itemize}
\item[(ii)] Réciproquement, lorsque les conditions ci-dessus sont réunies, la 2-bira\-tionalité de $L/K$ se propage à $L'/K'$.
\end{itemize}
\end{Th}

En particulier:

\begin{Sco}
La propagation de la 2-birationalité par 2-extension du corps de base ne peut se faire que par extension quadratique.
\end{Sco}

\begin{Sco}
Cette propagation ne peut se faire que si l'extension de départ $L/K$ est ramifiée modérément en deux places primitives.
\end{Sco}

\PreuveTh Supposons d'abord $L'/K'$  2-birationnelle.  Le corps $K'$ est alors une extension 2-rationnelle de $K$. Et, comme nous l'avons supposée linéairement disjointe de la 2-extension cyclotomique $K^c/K$, le Théorème \ref{Propa2r} nous assure que l'extension $K'/K$ est modérément ramifiée en une unique place modérée $\p$, laquelle est primitive dans $K$. Plus précisément, il résulte de l'hypothèse de disjonction linéaire que la place $\p$ est totalement ramifiée dans $K'/K$ (puisque la pro-2-extension 2-ramifiée $\infty$-décomposée maximale de $K$ est $K^c$). En particulier, l'unique place $\p'$ de $K'$ au-dessus de $\p$ est primitive dans $K'$.\par
Considérons le schéma d'extensions:
\begin{displaymath}
 \xymatrix { 
    &L=K[\sqrt{-\delta}] \ar@{-}[r] \ar@{-}[d] & L'=LK'\ar@{-}[d] \\
    &K \ar@{-}[r] & K'\\
    }
\end{displaymath}
Et examinons les deux possibilités décrites dans le Théorème \ref{Carac2bi} $(ii)$:\medskip

$\bullet$ $L'/K'$ est ramifiée modérément en une unique place semi-primitive $\q'$.\smallskip

\noindent Dans ce cas, $\q'$, qui n'est pas primitive, est distincte de $\p'$ et provient d'une place $\q$ de $K$ qui se ramifie dans $L/K$ mais ne se décompose pas dans $K'/K$. Ainsi $\q$ est (totalement) inerte dans $K'/K$, ce qui montre déjà que $K'/K$ est nécessairement cyclique. Plus précisément, puisque $\q'$ est semi-primitive dans $K'$, son degré d'inertie est (au plus) 2 et $\q$ doit être primitive. Ainsi d'une part $K'/K$ est quadratique; et d'autre part $L/K$, qui est ramifiée en la place primitive $\q$, dès lors qu'elle est 2-birationnelle se ramifie aussi en une autre place primitive (toujours en vertu du Théorème \ref{Carac2bi} (ii)), laquelle ne peut être que $\p$. 
\medskip

$\bullet$ $L'/K'$ est ramifiée modérément en deux places primitives $\q'_1$ et $\q'_2$.\smallskip

\noindent Dans ce cas, $\q'_1$ et $\q'_2$ proviennent soit d'une même place primitive $\q$ de $K$, soit de deux places primitives $\q_1$ et $\q_2$ distinctes, qui se ramifient dans $L/K$.\par
Dans cette dernière hypothèse, ni $\q_1$ ni $\q_2$ ne pourraient se décomposer dans $K'/K$, car $L'/K'$ serait alors ramifiée modérément en quatre places; elles ne pouraient non plus présenter de l'inertie, car $\q'_1$ ou $\q'_2$ ne serait plus primitive dans $K'$; elles seraient donc toutes deux totalement ramifiées dans $K'/K$, contrairement au fait que $\p$ est la seule place modérée ramifiée dans $K'/K$. \par
Il en résulte que $\q'_1$ et $\q'_2$ proviennent nécessairement d'une même place primitive $\q$ de $K$ dont l'indice de décomposition dans  $K'/K$ est égal à 2 et le degré d'inertie 1, puisque $\q'_1$ et $\q'_2$ sont primitives dans $K'$.
De plus, $\q$ est distincte de $\p$, sans quoi  $\q'_1$ et $\q'_2$ ne seraient pas ramifiées dans $L'/K'$, donc non ramifiée dans $K'/K$. \par
En résumé $[K':K]$ est égal à 2 et $\q'_1$ et $\q'_2$ proviennent d'une même place primitive de $K$ ramifiée dans $L/K$ et décomposée dans $K'/K$. Comme dans le cas précédent, $L/K$, qui est ramifiée en la place primitive $\q$, dès lors qu'elle est 2-birationnelle se ramifie nécessairement en une autre place primitive (toujours en vertu du Théorème \ref{Carac2bi} $(ii)$), laquelle ne peut être que $\p$. \medskip

Intéressons nous maintenant à la montée en supposant $L/K$ 2-birationnelle ramifiée modérément en deux places primitives $\p$ et $\q$, dont l'une, disons $\p$, est l'unique place modérée ramifiée dans l'extension quadratique $K'/K$. Et examinons successivement les deux possibilités recensées plus haut: \smallskip

$\bullet$ La place $\q$ est inerte dans $K'/K$.\smallskip

Dans ce cas, la place $\q'$ de $K'$ au-dessus de $\q$ est semi-primitive dans $K'$ et c'est l'unique place modérée ramifiée dans $L'/K'$, puisque $\p$, qui est totalement ramifiée dans $K'/K$, ne peut l'être aussi dans $L'/K$, son sous-groupe d'inertie étant cyclique. Ainsi $K'$ est bien 2-rationnel en vertu du Théorème \ref{Propa2r} et $L'$ est 2-birationnel d'après le Théorème \ref{Carac2bi} (ii).
\medskip

$\bullet$ La place $\q$ est décomposée dans $K'/K$.\smallskip

Dans ce cas, les deux place $\q'_1$ et $\q'_2$ de $L'$ au-dessus de $\q$ sont encore primitives et ce sont les seules places modérées qui se ramifient dans $L'/K'$, puisque, comme précédemment $\p$, qui est totalement ramifiée dans $K'/K$, ne peut l'être aussi dans $L'/K$. On conclut comme plus haut que $L'/K'$ est 2-birationnelle.

%%%%%%%%%%%%%%%%%%%%%%%%%%%%%%%%%%%%%%%%%%%%%%%%%%%%%%%

\section{\large Tours d'extensions 2-birationnelles}

Le Théorème \ref{ThPropa} limite aux seules extensions quadratiques la possibilité de réaliser la propagation de la 2-rationalité par 2-extension (galoisienne) du corps de base $K$. Mais il laisse ouverte la possibilité de construire des tours infinies (donc non galoisiennes) de telles extensions. Comme chaque montée quadratique n'est possible qu'à partir d'une extension birationnelle $L/K$ biramifiée modérément, la seule contrainte pour que la construction puisse se poursuivre est de ne fabriquer à chaque étage $n$ que des extension birationnelles $L_n/K_n$ qui soient encore biramifiées.\smallskip

En d'autres termes,  la question de la propagation indéfinie de la $2$-birationa\-lité se ramène à construire, pour une extension 2-birationnelle $L/K$ donnée, ramifiée modérément en exactement deux places primitives $\p$ et $\q$, une extension quadratique totalement réelle $K'$ de $K$ décomposée en $\q$ et ramifiée modérément en $\p$ seulement (la ramification sauvage étant indifférente).
\medskip

En vertu de la théorie 2-adique du corps de classes, une telle extension existe si et seulement si le $2$-groupe des $\infty\mathfrak{q}$-classes $2\mathfrak{p}$-infinitésimales est non trivial. \medskip

Regardons si c'est bien le cas. Rappelons le contexte de notre étude:\smallskip

\begin{enumerate}
\item $K$ est un corps de nombres totalement réel qui est $2$-rationnel; ce qui peut se traduire par les deux propriétés suivantes : 
\begin{enumerate}
     \item  $K$ possède une seule place dyadique $\l$; d'où: $[K_{\mathfrak{l}}:\Q_2]=[K:\Q]=r$;
     \item sa 2-extension abélienne 2-ramifiée $\infty$-décomposée maximale est sa $\Z_2$-extension cyclotomique $K^c$; ainsi le groupe d'idèles qui définit l'extension cyclotomique est donnée par la formule suivante: 
\begin{center}
                       $\wi\J_K=\prod\limits_{\mathfrak{r}\neq\mathfrak{l}}\mu_\mathfrak{r}\R_K$
\end{center}
\end{enumerate}
\item $L$ est une extension quadratique totalement imaginaire de $K$  qui se ramifie modérément en exactement deux places $\p$ et $\q$, qui sont primitives. En termes idéliques, cette primitivité s'écrit:
\begin{center}
$\J_K=\wi\J_K \R_{K\p}=\wi\J_K \R_{K_\q}$
\end{center}
\end{enumerate}\medskip

\noindent Cela étant, nous cherchons une extension quadratique $K'/K$ satisfaisant les quatre propriétés suivantes : 
\begin{enumerate}
\item[(i)] elle est ramifiée modérément en $\mathfrak{p}$ seulement,
\item[(ii)] elle est décomposée en $\mathfrak{q}$,
\item[(iii)] elle est non décomposée en $2$,
\item[(iv)] elle est complètement décomposée à l'infini (\ie totalement réelle).
\end{enumerate}\medskip

Pour l'instant, considérons la $2$-extension abélienne maximale $N$ de $K$ qui est totalement réelle, $\mathfrak{p}$-modérément ramifiée et $\mathfrak{q}$-décomposée. Le sous-groupe d'idèles qui lui correspond est ainsi: 
\begin{center}
$\prod\limits_{\r\mid\infty}\mu_\mathfrak{r}(\prod\limits_{\mathfrak{r}\nmid\mathfrak{l}\p\infty}\mu_\mathfrak{r})\mathcal{R}_{K_\q}\mathcal{R}_K$.
\end{center}       
Et il vient donc:
 \begin{center}
$\Gal(N/K)=\J_K/\prod\limits_{\r\nmid\p\l}\mu_\r\R_{K_\q}\R_K\simeq 
\mu_\p/ \mu_\p\cap(\prod\limits_{\r\nmid\p\l}\mu_\r \R_{K_\q}\R_K)$,
\end{center}
en vertu des égalités rappelées plus haut:
\begin{center}
$\J_K=\wi\J_K \R_{K_\q}$\qquad et \qquad $\wi\J_K =\prod\limits_{\r\neq\l}\mu_\r \R_K$.
\end{center}
Dans le quotient obtenu, les idèles principaux (\ie les éléments de $\mathcal{R}_K$) qui apparaissent au dénominateur sont dans $\prod\limits_{\mathfrak{r}\neq\mathfrak{l}}\mu_\mathfrak{r}$: ce sont des $\mathfrak{q}$-unités infinitésimales. Or, nous avons ici:

\begin{Lem}\label{ind} Dans un corps de nombres 2-rationnel totalement réel, pour toute place modérée $\q$ de $K$ le pro-2-groupe $\E^\q_\infty$ des $\q$-unités infinitésimales est trivial.
\end{Lem}

\Preuve Il s'agit de vérifier que l'image locale $s_\mathfrak l(\mathcal E^\mathfrak q )$ du 2-adifié $\mathcal E^\mathfrak q =\Z_2\otimes_\Z E^\mathfrak q$ du groupe des $\mathfrak q$-unités de $K$ est encore un $\Z_2$-module de rang $r=[K:\Q]$. Pour voir cela, observons que $K$,  puisqu'il est présumé 2-rationnel, vérifie la conjecture de Leopoldt; autrement dit que le 2-adifié groupe des unités $\mathcal E =\Z_2\otimes_\Z E$ s'injecte dans le groupe des unités locales $\mathcal U_\mathfrak l$ attaché à l'unique place dyadique $\mathfrak l$ de $K$. En particulier $\mathcal E$, qui est de rang $r-1 =[K:\Q]-1$, s'envoie avec un indice fini dans la préimage $\mathcal U_\mathfrak l^*$ dans $\mathcal U_\mathfrak l$ du groupe $\mu_2=\{\pm 1\}$ des racines de l'unité pour la norme arithmétique $\nu=N_{K/\Q}$. Soit alors $x$ l'image canonique dans $\mathcal E^\mathfrak q$ d'un générateur arbitraire d'une puissance principale de l'idéal $\mathfrak q$. La norme $x^\nu$ est (au signe près) une puissance non triviale de $N\mathfrak q$, et son logaritmhe 2-adique n'est donc pas nul. De sorte que le $\Z_2$-module $s_\l(\E^\q)$, qui contient $s_\l(\E)$ et $s_\l(x)$, est de rang au moins $(r-1)+1=r$; et finalement de rang exactement $r$, tout comme $\U_\l$.
En d'autres termes, le sous-module $\E^\q_\infty$ des $\q$-unités infinitésimales est bien trivial.\medskip

Ce point acquis, nous avons obtenu:

\begin{Prop}
Pour toute place primitive $\mathfrak q$ d'un corps 2-rationnel totalement réel $K$, la $2$-extension abélienne maximale $N$ de $K$ qui est totalement réelle, $\mathfrak{q}$-décomposée et ramifiée modérément en une unique place $\mathfrak{p}\ne \mathfrak q$, est cyclique, de groupe de Galois:

\centerline{$\Gal(N/K)\simeq\mu_\mathfrak p$.}\medskip

En particulier, $N/K$ contient une unique sous-extension quadratique $K'/K$.
\end{Prop}

Il suit de là que, sous les hypothèses de la proposition, l'unique sous-extension quadratique $K'/K$ de $N/K$ est l'unique extension quadratique qui satisfait aux conditions $(i)$, $(ii)$ et $(iv)$ listées plus haut. Reste à voir si elle vérifie également la condition $(iii)$ qui postule l'existence d'une unique place dyadique dans $K'$. Or c'est là qu'intervient précisément la condition de primitivité de la place $\mathfrak p$, que nous n'avons pas utilisée jusqu'ici: les résultats sur la propagation de la 2-rationalité rappelés dans le chapitre 1 assurent que $K'$ est encore 2-rationnel si et seulement si la place $\p$ est primitive dans $K$. Lorsque c'est le cas, $K'$ ne peut alors contenir qu'une seule place dyadique; et la condition $(iii)$ est, de ce fait, automatiquement vérifiée.\medskip

L'ensemble de cette discussion peut donc se résumer comme suit:

\begin{Th}\label {ThPrinc}
Soient $K$ un corps 2-rationnel totalement réel et $L$ une extension quadratique 2-birationnelle totalement imaginaire de $K$ ramifiée modérément en deux places primitives $\p$ et $\q$.
Il existe alors exactement deux extensions quadratiques $K'/K$ totalement réelles, 2-rationnelles et ramifiées modérément, telles que l'extension composée $L'=LK'$ soit 2-birationnelle:  celle qui est ramifiée modérément en $\p$ et décomposée en $\q$; et celle qui est  ramifiée modérément en $\q$ et décomposée en $\p$.
\end{Th}

Comme vu plus avant, l'extension $L'/K'$ vérifie à son tour les mêmes hypothèses que l'extension de départ. Itérant le théorème, on obtient ainsi:

\begin{Sco}
Sous les hypothèses du théorème, il existe une infinité de tours infinies d'extensions relativement quadratiques $K\subset K_1\subset K_2 \subset \cdots \subset K_i \subset \cdots$ de corps 2-rationnels totalement réels telles que les extensions composées $L_i=LK_i$ pour $i\in \N$ soient 2-birationnelles.
\end{Sco}

Il convient, en effet, à chaque étage $i\in \N$ déjà construit, de choisir celle des deux places primitives du corps $K_i$ ramifiées dans $L_i/K_i$ qu'on autorise à se ramifier modérément dans l'extension quadratique $K_{i+1}/K_i$.
\medskip

%%%%%%%%%%%%%%%%%%%%%%%%%%%%%%%%%%%%%%%%%%%%%%%%%%%%%%%%%
\section{\large Appendice: identités du miroir}
%%%%%%%%%%%%%%%%%%%%%%%%%%%%%%%%%%%%%%%%%%%%%%%%%%%%%%%%%

Il peut être instructif de relire les résultats ci-dessus à la lumière des identités du miroir qui fournissent une seconde preuve du Théorème \ref{ThPrinc} 
\medskip

Reprenons pour cela les calculs effectués plus haut: l'isomorphisme donné par le corps de classses

 \begin{center}
$\Gal(N/K)\simeq\mu_{\mathfrak{p}}/\mu_{\mathfrak{p}}\cap(\prod\limits_{\mathfrak{r}\nmid\mathfrak{p}\mathfrak{l}}\mu_\mathfrak{r}\mathcal{R}_\mathfrak{q}\mathcal{R}_K)$,
\end{center}

\noindent nous assure que l'extension $N/K$ est cyclique (éventuellement triviale). Posons $S=\lbrace\mathfrak{q}\infty\rbrace$ et $T=\lbrace\mathfrak{l}\mathfrak{p}\rbrace$. Par construction, le groupe de Galois $\Gal(N/K)$ s'identifie alors au 2-groupe des $S$-classes $T$-infinitésimales $\Cl_T^S$; et le résultat précédent se lit tout simplement:
$$
rg_2\,\Cl_T^S\leq 1.
$$

Nous allons à présent minorer ce rang grâce à la formule de réflexion de Gras (\cf Th. 4.6, p. 45 de \cite{Gr2}).

Reprenant les notations de Gras, nous avons: $S=\lbrace\mathfrak{q}\infty\rbrace$, $T=\lbrace\mathfrak{l}\mathfrak{p}\rbrace$, $S_0=\lbrace\mathfrak{q}\rbrace$, $T_2=\lbrace\mathfrak{l}\rbrace$, $S_2=\emptyset$ et $\Delta_\infty=\emptyset$; donc:
$$
rg_2\,\mathcal{C}\ell_T^S\;-\;rg_2\,\mathcal{C}\ell_{\lbrace\mathfrak{q}\rbrace}^{\lbrace\mathfrak{l}\mathfrak{p}\rbrace}= |T|+[K_\mathfrak l:\Q_2]-r-|S_0|-|\Delta_\infty|=2+r-r-1-0=1
$$
De cette formule, il suit en particulier: 
$$
rg_2\,\mathcal{C}\ell_T^S\geq1;
$$ 
de sorte qu'en fin de compte nous avons simultanément:
$$
rg_2\,\mathcal{C}\ell_T^S=1 \qquad {\rm et} \qquad \mathcal{C}\ell_{\lbrace\mathfrak{q}\rbrace}^{\lbrace\mathfrak{l}\mathfrak{p}\rbrace}=1.
$$ 

Le groupe $\mathcal{C}\ell_T^S$ est donc cyclique mais non trivial (comme nous l'avons déjà établi à l'aide du lemme d'indépendance \ref{ind} plus haut), tandis que le $\ell$-groupe $\mathcal{C}\ell_{\lbrace\mathfrak{q}\rbrace}^{\lbrace\mathfrak{l}\mathfrak{p}\rbrace}$ des $\{\mathfrak{lp}\}$-classes $\{\mathfrak q\}$-infinitésimales, lui, est nécessairement trivial.
\smallskip

De l'identité $rg_2\,\mathcal{C}\ell_T^S=1$, on conclut qu'il existe une unique extension quadratique $K'/K$ qui est non-ramifiée modérément en dehors de $\mathfrak{p}$ et $\infty\mathfrak{q}$-décomposée. Il reste alors à vérifier que cette extension est effectivement ramifiée en $\mathfrak{p}$ et qu'elle est non décomposée en $2$, pour qu'elle réalise la propagation de la $2$-birationalité.Or:\smallskip

\begin{itemize}
\item Le premier point est évident, puisque le groupe de Galois $\Gal(N/K) $ est engendré par l'image du sous-groupe d'inertie de la place $\mathfrak{p}$.

\item Il reste à voir que la place dyadique $\mathfrak l$ est non décomposée. Pour cela, reprenons le raisonnement précédent en échangeant les rôles de $\mathfrak{p}$ et de $\mathfrak{q}$.
Ce faisant, nous obtenons: $rg_2\,\mathcal{C}\ell_{\lbrace\mathfrak{p}\rbrace}^{\lbrace\mathfrak{l}\mathfrak{q}\rbrace}=0$, i.e. $\mathcal{C}\ell_{\lbrace\mathfrak{q}\rbrace}^{\lbrace\mathfrak{l}\mathfrak{p}\rbrace}=1$; ce qui est précisément le résultat attendu.
\end{itemize}\medskip

En conclusion, nous retrouvons le fait que pour un couple de places $\mathfrak{p}$ et $\mathfrak{q}$ primitives fixé, il existe une unique extension $K'/K$ vérifiant les hypothèses du théorème \ref{ThPropa} et permettant la propagation de la $2$-birationalité. \smallskip

De ce fait, ce processus peut être réitéré à l'infini et ainsi nous pouvons  construire des extensions $K'$ de $K$ de degrés arbitrairement grands, disjointes de la $\Z_2$-extension cyclotomique, de telle manière que le corps $L'$ obtenu soit encore $2$-birationnel.

\bigskip\bigskip

\noindent{\bf Adresses des auteurs}\medskip

\noindent Claire {\sc Bourbon},
Univ. Bordeaux,
Institut de Mathématiques de Bordeaux, UMR CNRS 5251,
351 Cours de la Libération,
F-33405 Talence cedex

\noindent{\small \tt claire.bourbon@math.u-bordeaux1.fr}
\medskip

\noindent Jean-François {\sc Jaulent},
Univ. Bordeaux,
Institut de Mathématiques de Bordeaux, UMR CNRS 5251,
351 Cours de la Libération,
F-33405 Talence cedex

\noindent{\small \tt jean-francois.jaulent@math.u-bordeaux1.fr}


\bigskip

\begin{thebibliography}{99}

\bibitem{Bo}{\sc C. Bourbon},
\textit{Propagation de la 2-rationalité},
{Thèse, Univ. Bordeaux, 2011}.

\bibitem{Gr1}{\sc G. Gras}, 
\textit{Théorèmes de réflexion},
{J. Théor. Nombres de Bordeaux} {\bf 10} (1998), 399--499.

\bibitem{Gr2}{\sc G. Gras}, 
\textit{Class Field Theory: From Theory To Practice},
{Springer-Verlag, 2003}.

\bibitem{GJ} {\sc G. Gras \& J.-F. Jaulent},
\textit{Sur les corps de nombres réguliers}, 
{Math. Z.}  {\bf 202} (1989), 343--365.

\bibitem{Ja1} {\sc J.-F Jaulent},
\textit{Classes logarithmiques des corps de nombres},
J. Théor. Nombres Bordeaux {\bf 10} (1994), 301--325.

\bibitem{Ja2}{\sc J.-F. Jaulent},
\textit {Théorie $\ell$-adique globale du corps de classes},
J. Théor. Nombres Bordeaux {\bf 10} (1998), 355--397.

\bibitem{JN} {\sc J.-F. Jaulent \& T. Nguyen Quang Do},
{\it Corps p-réguliers, corps p-rationnels et ramification restreinte}, 
J. Théor. Nombres Bordeaux  {\bf 5}  (1993), 343--363.

\bibitem{JS1} {\sc J.-F. Jaulent \& O. Sauzet},
\textit{Pro-$\ell$-extension de corps $\mathfrak l$-rationnels},
{J. Number Th.} {\bf 65} (1997), 240--267; {\it ibid.} {\bf 80} (2000), 318--319.

\bibitem{JS2}{\sc J.-F. Jaulent \& O. Sauzet},
\textit{Extensions quadratiques 2-birationnelles de corps totalement réels}, 
{Pub. Matemàtiques} {\bf 44} (2000), 343--351.

\bibitem{Mo} {\sc A. Movahhedi},
\textit{Sur les $p$-extensions des corps $p$-rationnels},
{Math. Nachr.} {\bf 149} (1990) 163--176.

\bibitem{MN} {\sc A. Movahhedi \& T. Nguyen Quang Do}, 
\textit{Sur l'arithmétique des corps de nombres $p$-rationnels},
{Sém. Th. Nombres Paris (1987/1988)}, Prog. in Math. {\bf 89} (1990), 155--200.

\bibitem {So}{\sc F. Soriano,} 
\textit{Classes logarithmiques ambiges des corps quadratiques}, 
{Acta Arith. } {\bf  78} (1997), 201--219.

\bibitem {Wi} {\sc  K. Wingberg,} 
{\it On Galois groups of p-closed algebraic number fields with restricted ramification},  
{J. reine  angew. Math.} {\bf 400} (1989), 185--202; {\it ibid.} {\bf 416} (1991),  187--194. 
 
\end{thebibliography}
\end{document}